\documentclass{article}
\usepackage{amssymb}
\usepackage{amsfonts}
\usepackage{amsmath}

\newtheorem{theorem}{Theorem}
\newtheorem{acknowledgement}[theorem]{Acknowledgement}

\newenvironment{proof}[1][Proof]{\noindent\textbf{#1.} }{\ \rule{0.5em}{0.5em}}
\input{tcilatex}

\begin{document}

\title{On the Lyapunov Exponent of Monotone Boolean Networks}
\author{Ilya Shmulevich \\
Institute for Systems Biology, Seattle, WA 98103\\
ilya.shmulevich@isbscience.org}

\maketitle

\begin{center}
    {\textit{Dedicated to the memory of A. D. Korshunov (1936-2019)}}
\end{center}

\begin{abstract}
We derive asymptotic formulas for the Lyapunov exponent of almost all monotone Boolean networks. The formulas are different depending on whether the number of variables of the constituent Boolean functions, or equivalently, the connectivity of the Boolean network, is even or odd.
\end{abstract}

\section{Introduction}

Boolean networks are complex dynamical systems that were proposed as models of
genetic regulatory networks \cite{kauffman1969a,kauffmanbook} and have
since then been used to model a range of complex phenomena. Random Boolean networks (RBNs) are ensembles of randomly generated Boolean networks with random topology and random updating functions. For a bias $p$ of the random Boolean functions, which is the probability that the function takes on the value 1, the critical connectivity is equal to $K_{c}=\left[  2p\left(
1-p\right)  \right]  ^{-1}$ (ref. \cite{derridapomeau}). Under a synchronous updating scheme, for $K<K_{c}$,
perturbations decay while for $K>K_{c}$, perturbations
spread throughout the network. This relationship is valid when the number of nodes in the network is infinite and coincides with the Lyapunov exponent \cite{luque,ilya-kauffman}. Networks that are known to operate in the dynamically critical regime, meaning that their Lyapunov exponent is $0$, are widely known to have many optimal properties and are thought to be a hallmark of living systems \cite{kauffmanbook,munoz}. 

We consider updating rules that belong to the class of monotone Boolean functions. This class of functions is one of the most widely studied classes of Boolean functions \cite{korshunov2003monotone}. Monotone Boolean networks have been primarily studied in the asynchronous updating scheme setting, whereby only one node is updated at a time. Some work has focused on long-term dynamics, such as fixed points and limit cycles \cite{aracena1,aracena2,aracena3}. It is also known that certain classes of fully asynchronous Boolean networks can be simulated by monotone Boolean networks \cite{melliti}. However, in the context of the Lyapunov exponent, asynchronous updating schemes appear to be less relevant than synchronous updating schemes \cite{mesot-teuscher}.

In the remainder, we will use $n$, rather than the conventional $K$ in the Boolean network literature, to denote the number of variables of the Boolean functions. Our results concerning almost all monotone Boolean networks can be understood in a probabilistic manner, meaning that the asymptotic formulas are valid with probability almost 1 if a monotone Boolean network is chosen at random from the set of all such networks.

As the formulas are asymptotic, they should not be interpreted for small $n$. At the same time, they are quite accurate even for $n=7$, which follows from the results originally published in \cite{korshunov}. Absent additional constraints on the Boolean functions, such as the classes of canalizing functions \cite {harris,serra,li,ilya-kauffman} or Post classes \cite{post}, random monotone Boolean networks quickly enter the disordered regime relative to $n$, but slower than RBNs \cite{derridapomeau}. Our main results here are asymptotic formulas, depending on whether $n$ is even or odd, for the so-called expected average sensitivity of a monotone Boolean function, first given in \cite{ilya-arxiv-MBF}. Given the results in \cite{ilya-kauffman}, the logarithm of the sensitivities, $\hat{s}^{f}$, in Theorems \ref{theven} and \ref{thodd}, can be directly interpreted as the Lyapunov exponent.

\section{Definitions and Preliminaries}

Let $f:\left\{ 0,1\right\} ^{n}\rightarrow \left\{ 0,1\right\} $ be a
Boolean function of $n$ variables $x_{1},\ldots ,x_{n}$. Let 
\begin{equation}
\partial f( \tilde{x}) /\partial x_{j}=f( \tilde{x}^{(
j,0) }) \oplus f( \tilde{x}^{( j,1) })
\end{equation}%
be the partial derivative of $f$ with respect to $x_{j}$, where $\oplus $ is
addition modulo 2 (exclusive OR) and $\tilde{x}^{\left( j,k\right) }=\left(
x_{1},\ldots ,x_{j-1},k,x_{j+1},\ldots x_{n}\right) $, $\,k=0,1$. Clearly,
the partial derivative is a Boolean function itself that specifies whether a
change in the $j$th input causes a change in the original function $f$. Now,
the activity of variable $x_{j}$ in function $f$ can be defined as 
\begin{equation}
\alpha _{j}^{f}=\frac{1}{2^{n}}\sum_{\tilde{x}\in \left\{ 0,1\right\}
^{n}}\partial f\left( \tilde{x}\right) /\partial x_{j}.  \label{eqactivity}
\end{equation}%
Note that although the vector $\tilde{x}$ consists of $n$ components
(variables), the $j$th variable is fictitious in $\partial f( \tilde{x}%
) /\partial x_{j}$. A variable $x_{j}$ is fictitious in $f$ if $%
f( \tilde{x}^{( j,0) }) =f( \tilde{x}^{(
j,1) }) $ for all $\tilde{x}^{( j,0) }$ and $\tilde{x}%
^{( j,1) }$. For a $n$-variable Boolean function $f$, we can form
its activity vector $\mathbf{\alpha }^{f}=[ \alpha _{1}^{f},\ldots
,\alpha _{n}^{f}] $. It is easy to see that $0\leq \alpha _{j}^{f}\leq
1,$ for any $j=1,\ldots ,n$. In fact, we can consider $\alpha _{j}^{f}$ to
be a probability that toggling the $j$th input bit changes the function
value, when the input vectors $\tilde{x}$ are distributed uniformly over $%
\{ 0,1\} ^{n}$. Since we're in the binary setting, the activity
is also the expectation of the partial derivative with respect to the
uniform distribution: $\alpha _{j}^{f}=E[ \partial f( \tilde{x}%
) /\partial x_{j}] $. Under an arbitrary distribution, $\alpha
_{j}^{f}$ is referred to as the \textit{influence} of variable $x_{j}$ on
the function $f$ \cite{kahn}. The influence of variables was used in the
context of genetic regulatory network modeling in \cite{pbn1}.

Another important quantity is the \textit{sensitivity} of a Boolean function 
$f$, which measures how sensitive the output of the function is to changes
in the inputs. The sensitivity $s^{f}\left( \tilde{x}\right) $ of $f$ on
vector $\tilde{x}$ is defined as the number of Hamming neighbors of $\tilde{x%
}$ on which the function value is different than on $\tilde{x}$ (two vectors
are Hamming neighbors if they differ in only one component). That is, 
\begin{align}
s^{f}\left( \tilde{x}\right) & =\left\vert \left\{ i\in \left\{ 1,\ldots
,n\right\} :f\left( \tilde{x}\oplus e_{i}\right) \neq f\left( \tilde{x}%
\right) \right\} \right\vert \\
& =\sum_{i=1}^{n}\chi \left[ f\left( \tilde{x}\oplus e_{i}\right) \neq
f\left( \tilde{x}\right) \right] ,  \notag
\end{align}%
where $e_{i}$ is the unit vector with 1 in the $i$th position and 0s
everywhere else and $\chi \left[ A\right] $ is an indicator function that is
equal to 1 if and only if $A$ is true. The \textit{average sensitivity} $%
s^{f}$ is defined by taking the expectation of $s^{f}\left( \tilde{x}\right) 
$ with respect to the distribution of $\tilde{x}$. It is easy to see that
under the uniform distribution, the average sensitivity is equal to the sum
of the activities: 
\begin{align}
s^{f}& =E\left[ s^{f}\left( \tilde{x}\right) \right] =\sum_{i=1}^{n}E\left[
\chi \left[ f\left( \tilde{x}\oplus e_{i}\right) \neq f\left( \tilde{x}%
\right) \right] \right]  \label{eqavgsens} \\
& =\sum_{i=1}^{n}\alpha _{i}^{f}.  \notag
\end{align}%
Therefore, $s^{f}$ is a number between $0$ and $n$.

Let $\tilde{\alpha}=\left( \alpha _{1},\cdots ,\alpha _{n}\right) $ and $%
\tilde{\beta}=\left( \beta _{1},\cdots ,\beta _{n}\right) $ be two different 
$n$-element binary vectors. We say that $\tilde{\alpha}$ precedes $\tilde{%
\beta}$, denoted as $\tilde{\alpha}\prec \tilde{\beta}$, if $\alpha _{i}\leq
\beta _{i}$ for every $i$, $1\leq i\leq n$. If $\tilde{\alpha}\nprec \tilde{%
\beta}$ and $\tilde{\beta}\nprec \tilde{\alpha}$, then $\tilde{\alpha}$ and $%
\tilde{\beta}$ are said to be incomparable. Relative to the predicate $\prec 
$, the set of all binary vectors of a given length is a partially ordered
set. A Boolean function $f\left( x_{1},\cdots ,x_{n}\right) $ is called 
\textit{monotone} if for any two vectors $\tilde{\alpha}$ and $\tilde{\beta}$
such that $\tilde{\alpha}\prec \tilde{\beta}$, we have $f\left( \tilde{\alpha%
}\right) \leq f(\tilde{\beta})$. 

We denote by $M\left( n\right) $ the set of all monotone Boolean functions
of $n$ variables. Let $E^{n}$ denote the Boolean $n$-cube, that is, a graph
with $2^{n}$ vertices each of which is labeled by an $n$-element binary
vector. Two vertices $\tilde{\alpha}=\left( \alpha _{1},\cdots ,\alpha
_{n}\right) $ and $\tilde{\beta}=\left( \beta _{1},\cdots ,\beta _{n}\right) 
$ are connected by an edge if and only if the Hamming distance $\rho (\tilde{%
\alpha},\tilde{\beta})=\sum_{i=1}^{n}\left( \alpha _{i}\oplus \beta
_{i}\right) =1$. The set of those vectors from $E^{n}$ in which there are
exactly $k$ units, $0\leq k\leq n$, is called the $k$th \textit{layer} of $%
E^{n}$ and is denoted by $E^{n,k}$.

A vector $\tilde{\alpha}\in E^{n}$ is called a \textit{minimal one} or \textit{minimal unit} of
monotone Boolean function $f( x_{1},\ldots ,x_{n}) $ if $f( 
\tilde{\alpha}) =1$ and $f( \tilde{\beta}) =0$ for any $%
\tilde{\beta}\prec \tilde{\alpha}$. A vector $\tilde{\alpha}\in E^{n}$ is
called an \textit{maximal zero}\emph{\ }of monotone Boolean function $%
f( x_{1},\ldots ,x_{n}) $ if $f( \tilde{\alpha}) =0$
and $f( \tilde{\beta}) =1$ for any $\tilde{\beta}\succ \tilde{%
\alpha}$. The minimal ones correspond directly to the terms in the minimal
disjunctive normal form (DNF) representation of the monotone Boolean
function. In \cite{ilya-normal}, asymptotic formulae for the number of
monotone Boolean functions of $n$ variables with a most probable number of
minimal ones were derived, confirming the conjecture in \cite{ilyansip95} that the number of monotone Boolean functions
relative to the number of minimal ones asymptotically follows a normal
distribution, with the assumption of all monotone Boolean functions being
equiprobable.

The average sensitivity has been studied intensively by a number of authors 
\cite{friedgutkalai,bernasconi,bernasconi2,boppana,friedgut,ilya-kauffman,shi,mossel,zhang,bshouty}. For example, it was shown by Friedgut 
\cite{friedgut} that if the average sensitivity of $f$ is $k$ then $f$ can
be approximated by a function depending on only $c^{k}$ variables where $c$
is a constant depending only on the accuracy of the approximation but not on 
$n$. Shmulevich and Kauffman \cite{ilya-kauffman} have shown that the
average sensitivity determines the critical phase transition curve in random
Boolean networks. Shi \cite{shi} showed that the average sensitivity can
serve as a lower bound of quantum query complexity. Average sensitivity was
used to characterize the noise sensitivity of monotone Boolean functions by
Mossel and O'Donnell \cite{mossel}. Zhang \cite{zhang} gives lower and upper bounds of the average sensitivity of a monotone Boolean function. The upper bound is asymptotic to $\sqrt{n}$, which has been shown by Bshouty and Tamon \cite{bshouty}. Our main results are given in Theorems \ref%
{theven} and \ref{thodd}.

\subsection{The structure of special monotone Boolean functions}

We now briefly review some known results concerning the structure of
so-called \textit{special }monotone Boolean functions. Let $M_{0}\left(
n\right) $ denote the set of functions in $M\left( n\right) $ possessing the
following properties. If $n$ is even, then $M_{0}\left( n\right) $ contains
only functions $f\in M\left( n\right) $ such that all minimal ones of $f$
are situated in $E^{n,n/2-1}$, $E^{n,n/2}$, and $E^{n,n/2+1}$ while function 
$f$ is equal to $1$ on all vectors in $E^{n,n/2+2},\cdots ,E^{n,n}$. For odd 
$n$, $M_{0}\left( n\right) $ contains only functions $f\in M\left( n\right) $
such that all minimal ones of $f$ are situated in either $E^{n,\left(
n-3\right) /2}$, $E^{n,\left( n-1\right) /2}$, and $E^{n,\left( n+1\right)
/2}$ or $E^{n,\left( n-1\right) /2}$, $E^{n,\left( n+1\right) /2}$, and $%
E^{n,\left( n+3\right) /2}$. In the first case, $f\left( \tilde{\alpha}%
\right) =1$ for all $\tilde{\alpha}$ in $E^{n,\left( n+3\right) /2},\cdots
,E^{n,n}$ while in the second case, $f\left( \tilde{\alpha}\right) =1$ for
all $\tilde{\alpha}$ in $E^{n,\left( n+5\right) /2},\cdots ,E^{n,n}$.

Then, as shown in \cite{korshunov}, 
\begin{equation}
\lim_{n\rightarrow \infty }\frac{\left\vert M_{0}\left( n\right) \right\vert 
}{\left\vert M\left( n\right) \right\vert }=1,  \label{eq23}
\end{equation}%
which we denote by $\left\vert M_{0}\left( n\right) \right\vert \sim
\left\vert M\left( n\right) \right\vert .$ In \cite{ilya-special},
asymptotic formulae for the number of special functions from $M_{0}\left(
n\right) $ were established and subsequently used to characterize
statistical properties of a popular class of nonlinear digital filters
called stack filters \cite{ilya-almostall}. The set of these special
functions is denoted by $M_{0}^{1}\left( n\right) $ and, depending on
whether $n$ is even or odd, is defined differently. While we shall omit the
rather lengthy definitions of special functions, the result from \cite%
{ilya-special} that will be important to us is that $\left\vert
M_{0}^{1}\left( n\right) \right\vert \sim \left\vert M\left( n\right)
\right\vert $. In other words, almost all monotone Boolean functions are
special. We shall also need the following results.

Let us start with the case of even $n$. Let 
\begin{align}
r_{0} & =r_{0}\left( n\right) =v_{0}=v_{0}\left( n\right) =\left\lfloor 
\binom{n}{n/2-1}2^{-n/2-1}\right\rfloor ,  \label{eq24} \\
z_{0} & =\left\lfloor \frac{1}{2}\binom{n}{n/2}\right\rfloor .  \notag
\end{align}
Let $M_{0}^{1}\left( n,r,z,v\right) $ denote the set of functions $f\in
M_{0}^{1}\left( n\right) $ such that $f$ has $r$ minimal ones in $%
E^{n,n/2-1} $, $v$ maximal zeros in $E^{n,n/2+1}$, and $f$ is equal to $1$
on $z$ vertices in $E^{n,n/2}$. In \cite{ilya-special}, the following result
was proved.

\begin{theorem}
\label{th1}Let $n$ be even, 
\begin{equation}
r=r_{0}+k,\,z=z_{0}+u,\,v=v_{0}+t,  \label{eq25}
\end{equation}
where $r_{0},z_{0},v_{0}$ are defined in (\ref{eq24}). Then, for any $k$, $t$%
, and $u$ such that $\left| k\right| \leq n2^{n/4}$, $\left| t\right| \leq
n2^{n/4}$, $\left| u\right| \leq n2^{n/2}$, 
\begin{multline*}
\left| M_{0}^{1}\left( n,r,z,v\right) \right| \sim\sqrt{\frac{2^{n+1}}{%
\pi^{3}\binom{n}{n/2}^{3}}}\left| M\left( n\right) \right| \\
\times\exp\left\{ -\frac{2^{n/2}}{\binom{n}{n/2-1}}\left( k^{2}+t^{2}\right)
-\frac{2u^{2}}{\binom{n}{n/2}}\right\} .
\end{multline*}
\end{theorem}

For any odd $n$, we use the parameters $r_{1},z_{1},v_{1}$ which are given
by 
\begin{align}
r_{1} & =r_{1}\left( n\right) =\left\lfloor \binom{n}{\left( n-3\right) /2}%
2^{-\left( n+3\right) /2}\right\rfloor ,\,  \label{eq27} \\
v_{1} & =v_{1}\left( n\right) =\left\lfloor \binom{n}{\left( n+1\right) /2}%
2^{-\left( n+1\right) /2}\right\rfloor ,  \notag \\
z_{1} & =\left\lfloor \frac{1}{2}\left( \binom{n}{\left( n-1\right) /2}%
+r_{1}\left( n+3\right) /2-v_{1}\left( n+1\right) /2\right) \right\rfloor
\label{eq28}
\end{align}
and parameters $r_{2},z_{2},v_{2}$, which are given by 
\begin{align}
r_{2} & =r_{2}\left( n\right) =\left\lfloor \binom{n}{\left( n-1\right) /2}%
2^{-\left( n+1\right) /2}\right\rfloor ,\,  \label{eq210} \\
v_{2} & =v_{2}\left( n\right) =\left\lfloor \binom{n}{\left( n+3\right) /2}%
2^{-\left( n+3\right) /2}\right\rfloor ,  \notag \\
z_{2} & =\left\lfloor \frac{1}{2}\left( \binom{n}{\left( n+1\right) /2}%
+r_{2}\left( n-1\right) /2-v_{2}\left( n+3\right) /2\right) \right\rfloor .
\label{eq211}
\end{align}

Let $M_{0,1}^{1}\left( n,r,z,v\right) $ denote the set of functions $f\in
M_{0}^{1}\left( n\right) $ such that $f$ has $r$ minimal ones in $%
E^{n,\left( n-3\right) /2}$, $v$ maximal zeros in $E^{n,\left( n+1\right)
/2} $, and $f$ is equal to $1$ on $z$ vertices in $E^{n,\left( n-1\right)
/2} $. Similarly, let $M_{0,2}^{1}\left( n,r,z,v\right) $ denote the set of
functions $f\in M_{0}^{1}\left( n\right) $ such that $f$ has $r$ minimal
ones in $E^{n,\left( n-1\right) /2}$, $v$ maximal zeros in $E^{n,\left(
n+3\right) /2}$, and $f$ is equal to $1$ on $z$ vertices in $E^{n,\left(
n+1\right) /2}$. Then, in \cite{ilya-special}, the following two Theorems
were proved.

\begin{theorem}
\label{th2}Let $n$ be odd, 
\begin{equation}
r=r_{1}+k,\,z=z_{1}+u,\,v=v_{1}+t,  \label{eq29}
\end{equation}
where $r_{1},z_{1},v_{1}$ are defined in (\ref{eq27}) and (\ref{eq28}).
Then, for any $k$, $t$, and $u$ such that $\left| k\right| \leq n2^{n/4}$, $%
\left| t\right| \leq n2^{n/4}$, $\left| u\right| \leq n2^{n/2}$, 
\begin{multline*}
\left| M_{0,1}^{1}\left( n,r,z,v\right) \right| \sim\frac{1}{2}\sqrt {\frac{%
2^{n+1}}{\pi^{3}\binom{n}{\left( n-1\right) /2}^{3}}}\left| M\left( n\right)
\right| \\
\times\exp\left\{ -\frac{2^{\left( n+1\right) /2}}{\binom{n}{\left(
n-3\right) /2}}k^{2}-\frac{2^{\left( n-1\right) /2}}{\binom{n}{\left(
n+1\right) /2}}t^{2}-\frac{2u^{2}}{\binom{n}{\left( n-1\right) /2}}\right\} .
\end{multline*}
\end{theorem}

\begin{theorem}
\label{th3}Let $n$ be odd, 
\begin{equation}
r=r_{2}+k,\,z=z_{2}+u,\,v=v_{2}+t,  \label{eq212}
\end{equation}%
where $r_{2},z_{2},v_{2}$ are defined in (\ref{eq210}) and (\ref{eq211}).
Then, for any $k$, $t$, and $u$ such that $\left\vert k\right\vert \leq
n2^{n/4}$, $\left\vert t\right\vert \leq n2^{n/4}$, $\left\vert u\right\vert
\leq n2^{n/2}$, 
\begin{multline*}
\left\vert M_{0,2}^{1}\left( n,r,z,v\right) \right\vert \sim \frac{1}{2}%
\sqrt{\frac{2^{n+1}}{\pi ^{3}\binom{n}{\left( n-1\right) /2}^{3}}}\left\vert
M\left( n\right) \right\vert \\
\times \exp \left\{ -\frac{2^{\left( n-1\right) /2}}{\binom{n}{\left(
n-1\right) /2}}k^{2}-\frac{2^{\left( n+1\right) /2}}{\binom{n}{\left(
n+3\right) /2}}t^{2}-\frac{2u^{2}}{\binom{n}{\left( n+1\right) /2}}\right\} .
\end{multline*}
\end{theorem}

\section{Main Results}

Since $\left\vert M_{0}\left( n\right) \right\vert \sim \left\vert M\left(
n\right) \right\vert ,$ we can focus our attention on functions in $%
M_{0}\left( n\right) $ and derive the average sensitivity of a typical
function from $M_{0}\left( n\right) .$ By `typical' we mean the most
probable Boolean function relative to the parameters $k$, $t$, and $u$ in
Theorems \ref{th1}-\ref{th3}. It can easily be seen that the most probable
special Boolean functions will have $k=t=u=0.$ This will imply, to take the $%
n$-even case as an example, that the most probable function $f$ has $r_{0}$
minimal ones in $E^{n,n/2-1}$, $v_{0}$ maximal zeros in $E^{n,n/2+1}$, and $%
f $ is equal to $1$ on $z_{0}$ vertices in $E^{n,n/2}$, where $r_{0},v_{0},$
and $z_{0}$ are given in (\ref{eq24}). Our proofs are thus based on the
derivation of the average sensitivity of such a function. Whenever we make
probabilistic assertions using words such as `most probable' or `typical' or
talk about expectations, we are implicitly endowing the set $M_{0}^{1}\left(
n,r,z,v\right) $ with a uniform probability distribution for fixed
parameters $n,r,z,v.$ This should not be confused with the Gaussian-like
distribution of $M_{0}^{1}\left( n,r,z,v\right) $ relative to its parameters 
$n,r,z,v.$ We will also omit the floor notation $\left\lfloor \cdot
\right\rfloor $ as the results are asymptotic.

\begin{theorem}
\label{theven}Let $n$ be even and let $f\in M_{0}^{1}\left( n\right) $ be a
typical monotone Boolean function. Then, the expected average sensitivity $%
\hat{s}^{f}=E\left[ s^{f}\right] $ of $f$ is%
\begin{equation*}
\hat{s}^{f}\sim n2^{-n}\binom{n}{n/2-1}\left( 2^{-n/2-1}+1\right) .
\end{equation*}
\end{theorem}

\begin{proof}
We will proceed by first focusing on determining the activity of an
arbitrary variable $x_{j}$ of a typical function $f.$ By simple symmetry
arguments, if we were to sample randomly from the set $M\left( n\right) $ of
monotone Boolean functions, the expected activities would be equal for all
the variables. It will follow by (\ref{eqavgsens}) that the expected average
sensitivity will be equal to $n$ multiplied by the expected activity. Since
the function $f$ is such that its minimal ones are situated in $E^{n,n/2-1}$%
, $E^{n,n/2}$, and $E^{n,n/2+1}$ while it is equal to $1$ on all vectors in $%
E^{n,n/2+2},\cdots ,E^{n,n}$, the only non-trivial behavior occurs between
the layers $E^{n,n/2-2}$ and $E^{n,n/2+2}.$

Let us consider the minimal ones, and hence all of the ones, on $%
E^{n,n/2-1}. $ Since we are considering variable $x_{j},$ half of these
minimal ones will have $x_{j}=0$ (i.e. $\tilde{x}^{\left( j,0\right) }$) and
the other half will have $x_{j}=1$ (i.e. $\tilde{x}^{\left( j,1\right) }$)$.$
It is easy to see that if $\tilde{x}^{\left( j,0\right) }\in E^{n,n/2-1}$ is
a minimal one, then by monotonicity, $f\left( \tilde{x}^{\left( j,1\right)
}\right) =1$. Consequently, the Hamming neighbors $\tilde{x}^{\left(
j,0\right) }$ and $\tilde{x}^{\left( j,1\right) }$ contribute nothing to the
sum in (\ref{eqactivity}). On the other hand, if $\tilde{x}^{\left(
j,1\right) }\in E^{n,n/2-1}$ is a minimal one, then $\partial f\left( \tilde{%
x}\right) /\partial x_{j}=1,$ since $f\left( \tilde{x}\right) =0$ for all $%
\tilde{x}\in E^{n,n/2-2}$. The (most probable\footnote{%
To avoid repetition, we will omit the words "most probable" or "typical"
when it is understood from the context.}) number of such minimal ones on $%
E^{n,n/2-1}$ contributing to the sum in (\ref{eqactivity}) is thus equal to 
\begin{equation}
\frac{1}{2}\binom{n}{n/2-1}2^{-n/2-1}.  \label{eqA}
\end{equation}

The number of zeros on $E^{n,n/2-1}$ is equal to 
\begin{equation}
\binom{n}{n/2-1}-r_{0}=\binom{n}{n/2-1}\left( 1-2^{-n/2-1}\right) .
\label{eqnumzeros}
\end{equation}%
As above, half of these will have $x_{j}=0$ and half will have $x_{j}=1$. We
need not consider vectors $\tilde{x}^{\left( j,1\right) }\in E^{n,n/2-1},$
since $f\left( \tilde{x}\right) =0$ for all $\tilde{x}\in E^{n,n/2-2}.$
However, we should consider the number of ones situated on the middle layer $%
E^{n,n/2}.$ The middle layer contains 
\begin{equation}
\frac{1}{2}\binom{n}{n/2}
\end{equation}%
ones and an equal number of zeros. Thus, half of the vectors $\tilde{x}%
^{\left( j,1\right) }\in E^{n,n/2}$ will be ones and the other half will be
zeros. In total, the number of vectors $\tilde{x}\in E^{n,n/2-1}$ such that $%
f\left( \tilde{x}\right) =0,$ $x_{j}=0,$ and $f\left( \tilde{x}^{\left(
j,1\right) }\right) =1,$ is equal to 
\begin{equation}
\frac{1}{4}\binom{n}{n/2-1}\left( 1-2^{-n/2-1}\right) .  \label{eqB}
\end{equation}%
We have now examined all partial derivatives above and below the layer $%
E^{n,n/2-1}.$

Let us now jump to layer $E^{n,n/2+1},$ as it will be similar by duality
considerations. The number of maximal zeros on that layer is equal to $%
v_{0}=r_{0}.$ Half of these will have $x_{j}=1$ and thus $\partial f\left( 
\tilde{x}\right) /\partial x_{j}=0$ due to monotonicity. The other half will
have $x_{j}=0$ and since $f\left( \tilde{x}\right) =1$ for all $\tilde{x}\in
E^{n,n/2+2},$ the same total as in (\ref{eqA}) will result. Similarly, the
number of ones on $E^{n,n/2+1}$ is the same as in (\ref{eqnumzeros}). We are
only concerned with $\tilde{x}^{\left( j,1\right) }\in E^{n,n/2+1},$ such
that $f\left( \tilde{x}^{\left( j,1\right) }\right) =1$ and $f\left( \tilde{x%
}^{\left( j,0\right) }\right) =0.$ As above, because the middle layer
contains the same number of ones and zeros, the total number of such pairs
of vectors is the same as in (\ref{eqB}). Thus, having accounted for all
partial derivatives above and below $E^{n,n/2+1}$ and having convinced
ourselves that their contribution to the overall activity of variable $x_{j}$
is the same as in (\ref{eqA}) and (\ref{eqB}), we can multiply (\ref{eqA})
and (\ref{eqB}) by $2$ and add them together to obtain 
\begin{eqnarray}
&&\binom{n}{n/2-1}2^{-n/2-1}+\frac{1}{2}\binom{n}{n/2-1}\left(
1-2^{-n/2-1}\right) \\
&=&\frac{1}{2}\binom{n}{n/2-1}\left( 2^{-n/2-1}+1\right)  \label{eqC}
\end{eqnarray}
Since there are $2^{-n+1}$ Hamming neighbors $\tilde{x}^{\left( j,0\right)
}\prec \tilde{x}^{\left( j,1\right) }$ and since the average sensitivity is $%
n$ times the activity, we must multiply (\ref{eqC}) by $n2^{-n+1},$
resulting in the statement of the theorem.
\end{proof}

The case of odd $n$ is somewhat more involved because there is no "middle"
layer $E^{n,n/2}.$ Instead, typical functions break up into two sets: $%
M_{0,1}^{1}\left( n,r,z,v\right) $ and $M_{0,2}^{1}\left( n,r,z,v\right) .$
In the first case, all minimal ones are on layers $E^{n,\left( n-3\right)
/2} $, $E^{n,\left( n-1\right) /2}$, and $E^{n,\left( n+1\right) /2}$ while
in the second case, all minimal ones are situated on $E^{n,\left( n-1\right)
/2} $, $E^{n,\left( n+1\right) /2}$, and $E^{n,\left( n+3\right) /2}$. Under
random sampling, these two cases will occur with equal probabilities. As we
shall see, the results will be different for each case. Thus, the expected
average sensitivity will be the average of the the expected average
sensitivities corresponding to these two cases.

\begin{theorem}
\label{thodd}Let $n$ be odd and let $f\in M_{0}^{1}\left( n\right) $ be a
typical monotone Boolean function. Then, the expected average sensitivity $%
\hat{s}^{f}=E\left[ s^{f}\right] $ of $f$ is%
\begin{equation*}
\hat{s}^{f}\sim \frac{1}{2}\left( \hat{s}_{1}^{f}+\hat{s}_{2}^{f}\right) ,
\end{equation*}%
where $\hat{s}_{1}^{f}$ and $\hat{s}_{2}^{f}$ are given in (\ref{eqsf1}) and
(\ref{eqsf2}), respectively.
\end{theorem}

\begin{proof}
Let us first address the set $M_{0,1}^{1}\left( n,r,z,v\right) .$ As in
Theorem \ref{theven}, we only need to concern ourselves with four cases.
These are:

\begin{enumerate}
\item $\tilde{x}\in E^{n,\left( n-3\right) /2},$ $f\left( \tilde{x}\right)
=1,$ $x_{j}=1$

\item $\tilde{x}\in E^{n,\left( n-3\right) /2},$ $f\left( \tilde{x}\right)
=0,$ $x_{j}=0$

\item $\tilde{x}\in E^{n,\left( n+1\right) /2},$ $f\left( \tilde{x}\right)
=1,$ $x_{j}=1$

\item $\tilde{x}\in E^{n,\left( n+1\right) /2},$ $f\left( \tilde{x}\right)
=0,$ $x_{j}=0$
\end{enumerate}

All other situations, such as $\tilde{x}\in E^{n,\left( n-3\right) /2},$ $%
f\left( \tilde{x}\right) =1,$ $x_{j}=0,$ will result in the partial
derivatives being zero due to monotonicity of $f$, hence will make no
contribution to the activity of variable $x_{j}.$ There are 
\begin{equation}
\frac{1}{2}\binom{n}{\left( n-3\right) /2}2^{-\left( n+3\right) /2}
\label{eqD}
\end{equation}%
minimal ones on $E^{n,\left( n-3\right) /2}$ such that $x_{j}=1.$ At the
same time, there are 
\begin{equation}
\frac{1}{2}\binom{n}{\left( n-3\right) /2}\left( 1-2^{-\left( n+3\right)
/2}\right)  \label{eqE}
\end{equation}%
zeros on $E^{n,\left( n-3\right) /2}$ such that $x_{j}=0.$ Unlike in the $n$%
-even case, where the middle layer contains an equal number of ones and
zeros, the number of ones and zeros on $E^{n,\left( n-1\right) /2}$ is not
equal. The number of ones on $E^{n,\left( n-1\right) /2}$ is given in (\ref%
{eq28}). Thus, if $\tilde{x}^{\left( j,0\right) }$ is a zero on $E^{n,\left(
n-3\right) /2}$, then the probability that $f\left( \tilde{x}^{\left(
j,1\right) }\right) =1$ is 
\begin{multline}
\binom{n}{\left( n-1\right) /2}^{-1}\left( \frac{1}{2}\left( \binom{n}{%
\left( n-1\right) /2}+\binom{n}{\left( n-3\right) /2}2^{-\left( n+3\right)
/2}\left( n+3\right) /2\right. \right.  \label{eqF} \\
\left. \left. -\binom{n}{\left( n+1\right) /2}2^{-\left( n+1\right)
/2}\left( n+1\right) /2\right) \right)
\end{multline}%
where we have simply divided the number of ones on $E^{n,\left( n-1\right)
/2}$ by the total number of vectors on that layer. Thus, multiplying (\ref%
{eqE}) by (\ref{eqF}) and adding to (\ref{eqD}) gives us the total
contribution to the activity of variable $x_{j}$ from cases 1 and 2 above,
which is equal to%
\begin{multline}
\frac{1}{2}\binom{n}{\left( n-3\right) /2}2^{-\left( n+3\right) /2}+\frac{1}{%
2}\binom{n}{\left( n-3\right) /2}\left( 1-2^{-\left( n+3\right) /2}\right)
\times  \label{eqF2} \\
\binom{n}{\left( n-1\right) /2}^{-1}\left( \frac{1}{2}\left( \binom{n}{%
\left( n-1\right) /2}+\binom{n}{\left( n-3\right) /2}2^{-\left( n+3\right)
/2}\left( n+3\right) /2\right. \right. \\
\left. \left. -\binom{n}{\left( n+1\right) /2}2^{-\left( n+1\right)
/2}\left( n+1\right) /2\right) \right)
\end{multline}

Cases 3 and 4 are similar. The number of (maximal) zeros on $E^{n,\left(
n+1\right) /2}$ for which $x_{j}=0$ is equal to 
\begin{equation}
\frac{1}{2}\binom{n}{\left( n+1\right) /2}2^{-\left( n+1\right) /2}
\label{eqG}
\end{equation}%
and the number of ones on $E^{n,\left( n+1\right) /2}$ for which $x_{j}=1$
is equal to 
\begin{equation}
\frac{1}{2}\binom{n}{\left( n+1\right) /2}\left( 1-2^{-\left( n+1\right)
/2}\right) .  \label{eqH}
\end{equation}%
The proportion of zeros on $E^{n,\left( n-1\right) /2}$ is 
\begin{multline}
\left( 1-\binom{n}{\left( n-1\right) /2}^{-1}\left( \frac{1}{2}\left( \binom{%
n}{\left( n-1\right) /2}+\right. \right. \right.  \label{eqI} \\
\left. \left. \left. \binom{n}{\left( n-3\right) /2}2^{-\left( n+3\right)
/2}\left( n+3\right) /2-\binom{n}{\left( n+1\right) /2}2^{-\left( n+1\right)
/2}\left( n+1\right) /2\right) \right) \right)
\end{multline}%
Thus, as before, multiplying (\ref{eqH}) by (\ref{eqI}) and adding to (\ref%
{eqG}), we get the total contribution to the activity from above and below
the layer $E^{n,\left( n+1\right) /2}$ (cases 3 and 4), resulting in 
\begin{multline}
\frac{1}{2}\binom{n}{\left( n+1\right) /2}2^{-\left( n+1\right) /2}+\frac{1}{%
2}\binom{n}{\left( n+1\right) /2}\left( 1-2^{-\left( n+1\right) /2}\right)
\times \\
\left( 1-\binom{n}{\left( n-1\right) /2}^{-1}\left( \frac{1}{2}\left( \binom{%
n}{\left( n-1\right) /2}+\right. \right. \right. \\
\left. \left. \left. \binom{n}{\left( n-3\right) /2}2^{-\left( n+3\right)
/2}\left( n+3\right) /2-\binom{n}{\left( n+1\right) /2}2^{-\left( n+1\right)
/2}\left( n+1\right) /2\right) \right) \right)  \label{eqI2}
\end{multline}

Finally, adding (\ref{eqF2}) and (\ref{eqI2}) and then multiplying by $%
n2^{-n+1}$ as in Theorem \ref{theven}, we obtain the expected average
sensitivity $\hat{s}_{1}^{f}$ of a typical function from $M_{0,1}^{1}\left(
n,r,z,v\right) $:%
\begin{multline}
\hat{s}_{1}^{f}\sim n2^{-n+1}\left( \frac{1}{2}\binom{n}{\left( n-3\right) /2%
}2^{-\left( n+3\right) /2}+\frac{1}{2}\binom{n}{\left( n-3\right) /2}\left(
1-2^{-\left( n+3\right) /2}\right) \times \right.  \label{eqsf1} \\
\binom{n}{\left( n-1\right) /2}^{-1}\left( \frac{1}{2}\left( \binom{n}{%
\left( n-1\right) /2}+\binom{n}{\left( n-3\right) /2}2^{-\left( n+3\right)
/2}\left( n+3\right) /2\right. \right. \\
\left. \left. -\binom{n}{\left( n+1\right) /2}2^{-\left( n+1\right)
/2}\left( n+1\right) /2\right) \right) + \\
\frac{1}{2}\binom{n}{\left( n+1\right) /2}2^{-\left( n+1\right) /2}+\frac{1}{%
2}\binom{n}{\left( n+1\right) /2}\left( 1-2^{-\left( n+1\right) /2}\right)
\times \\
\left( 1-\binom{n}{\left( n-1\right) /2}^{-1}\left( \frac{1}{2}\left( \binom{%
n}{\left( n-1\right) /2}+\right. \right. \right. \\
\left. \left. \left. \left. \binom{n}{\left( n-3\right) /2}2^{-\left(
n+3\right) /2}\left( n+3\right) /2-\binom{n}{\left( n+1\right) /2}2^{-\left(
n+1\right) /2}\left( n+1\right) /2\right) \right) \right) \right) .
\end{multline}%
We now proceed to derive the expected average sensitivity $\hat{s}_{2}^{f}$
of a typical function from $M_{0,2}^{1}\left( n,r,z,v\right) ,$ following
the same steps.

Again, we only need to concern ourselves with four cases:

\begin{enumerate}
\item $\tilde{x}\in E^{n,\left( n-1\right) /2},$ $f\left( \tilde{x}\right)
=1,$ $x_{j}=1$

\item $\tilde{x}\in E^{n,\left( n-1\right) /2},$ $f\left( \tilde{x}\right)
=0,$ $x_{j}=0$

\item $\tilde{x}\in E^{n,\left( n+3\right) /2},$ $f\left( \tilde{x}\right)
=1,$ $x_{j}=1$

\item $\tilde{x}\in E^{n,\left( n+3\right) /2},$ $f\left( \tilde{x}\right)
=0,$ $x_{j}=0$
\end{enumerate}

There are 
\begin{equation}
\frac{1}{2}\binom{n}{\left( n-1\right) /2}2^{-\left( n+1\right) /2}
\label{eqJ}
\end{equation}%
minimal ones on $E^{n,\left( n-1\right) /2}$ such that $x_{j}=1.$ At the
same time, there are 
\begin{equation}
\frac{1}{2}\binom{n}{\left( n-1\right) /2}\left( 1-2^{-\left( n+1\right)
/2}\right)  \label{eqK}
\end{equation}%
zeros on $E^{n,\left( n-1\right) /2}$ such that $x_{j}=0.$ As for $%
M_{0,1}^{1}\left( n,r,z,v\right) $, the number of ones and zeros on $%
E^{n,\left( n+1\right) /2}$ is not equal. The number of ones on $E^{n,\left(
n+1\right) /2}$ is given in (\ref{eq211}). Thus, if $\tilde{x}^{\left(
j,0\right) }$ is a zero on $E^{n,\left( n-1\right) /2}$, then the
probability that $f\left( \tilde{x}^{\left( j,1\right) }\right) =1$ is 
\begin{multline}
\binom{n}{\left( n+1\right) /2}^{-1}\left( \frac{1}{2}\left( \binom{n}{%
\left( n+1\right) /2}+\binom{n}{\left( n-1\right) /2}2^{-\left( n+1\right)
/2}\left( n-1\right) /2\right. \right.  \label{eqL} \\
\left. \left. -\binom{n}{\left( n+3\right) /2}2^{-\left( n+3\right)
/2}\left( n+3\right) /2\right) \right)
\end{multline}%
where we have simply divided the number of ones on $E^{n,\left( n+1\right)
/2}$ by the total number of vectors on that layer. Thus, multiplying (\ref%
{eqK}) by (\ref{eqL}) and adding to (\ref{eqJ}) gives us the total
contribution to the activity of variable $x_{j}$ from cases 1 and 2 above,
which is equal to%
\begin{multline}
\frac{1}{2}\binom{n}{\left( n-1\right) /2}2^{-\left( n+1\right) /2}+\frac{1}{%
2}\binom{n}{\left( n-1\right) /2}\left( 1-2^{-\left( n+1\right) /2}\right)
\times  \label{eqL2} \\
\binom{n}{\left( n+1\right) /2}^{-1}\left( \frac{1}{2}\left( \binom{n}{%
\left( n+1\right) /2}+\binom{n}{\left( n-1\right) /2}2^{-\left( n+1\right)
/2}\left( n-1\right) /2\right. \right. \\
\left. \left. -\binom{n}{\left( n+3\right) /2}2^{-\left( n+3\right)
/2}\left( n+3\right) /2\right) \right)
\end{multline}

Cases 3 and 4 are similar. The number of (maximal) zeros on $E^{n,\left(
n+3\right) /2}$ for which $x_{j}=0$ is equal to 
\begin{equation}
\frac{1}{2}\binom{n}{\left( n+3\right) /2}2^{-\left( n+3\right) /2}
\label{eqM}
\end{equation}%
and the number of ones on $E^{n,\left( n+3\right) /2}$ for which $x_{j}=1$
is equal to 
\begin{equation}
\frac{1}{2}\binom{n}{\left( n+3\right) /2}\left( 1-2^{-\left( n+3\right)
/2}\right) .  \label{eqN}
\end{equation}%
The proportion of zeros on $E^{n,\left( n+1\right) /2}$ is 
\begin{multline}
\left( 1-\binom{n}{\left( n+1\right) /2}^{-1}\left( \frac{1}{2}\left( \binom{%
n}{\left( n+1\right) /2}+\right. \right. \right.  \label{eqO} \\
\left. \left. \left. \binom{n}{\left( n-1\right) /2}2^{-\left( n+1\right)
/2}\left( n-1\right) /2-\binom{n}{\left( n+3\right) /2}2^{-\left( n+3\right)
/2}\left( n+3\right) /2\right) \right) \right)
\end{multline}%
Thus, as before, multiplying (\ref{eqN}) by (\ref{eqO}) and adding to (\ref%
{eqM}), we get the total contribution to the activity from above and below
the layer $E^{n,\left( n+3\right) /2}$ (cases 3 and 4), resulting in 
\begin{multline}
\frac{1}{2}\binom{n}{\left( n+3\right) /2}2^{-\left( n+3\right) /2}+\frac{1}{%
2}\binom{n}{\left( n+3\right) /2}\left( 1-2^{-\left( n+3\right) /2}\right)
\times \\
\left( 1-\binom{n}{\left( n+1\right) /2}^{-1}\left( \frac{1}{2}\left( \binom{%
n}{\left( n+1\right) /2}+\right. \right. \right. \\
\left. \left. \left. \binom{n}{\left( n-1\right) /2}2^{-\left( n+1\right)
/2}\left( n-1\right) /2-\binom{n}{\left( n+3\right) /2}2^{-\left( n+3\right)
/2}\left( n+3\right) /2\right) \right) \right)  \label{eqO2}
\end{multline}

Finally, adding (\ref{eqL2}) and (\ref{eqO2}) and then multiplying by $%
n2^{-n+1}$ as in Theorem \ref{theven}, we obtain the expected average
sensitivity $\hat{s}_{2}^{f}$ of a typical function from $M_{0,2}^{1}\left(
n,r,z,v\right) $:%
\begin{multline}
\hat{s}_{2}^{f}\sim n2^{-n+1}\left( \frac{1}{2}\binom{n}{\left( n-1\right) /2%
}2^{-\left( n+1\right) /2}+\frac{1}{2}\binom{n}{\left( n-1\right) /2}\left(
1-2^{-\left( n+1\right) /2}\right) \times \right.   \label{eqsf2} \\
\binom{n}{\left( n+1\right) /2}^{-1}\left( \frac{1}{2}\left( \binom{n}{%
\left( n+1\right) /2}+\binom{n}{\left( n-1\right) /2}2^{-\left( n+1\right)
/2}\left( n-1\right) /2\right. \right.  \\
\left. \left. -\binom{n}{\left( n+3\right) /2}2^{-\left( n+3\right)
/2}\left( n+3\right) /2\right) \right) + \\
\frac{1}{2}\binom{n}{\left( n+3\right) /2}2^{-\left( n+3\right) /2}+\frac{1}{%
2}\binom{n}{\left( n+3\right) /2}\left( 1-2^{-\left( n+3\right) /2}\right)
\times  \\
\left( 1-\binom{n}{\left( n+1\right) /2}^{-1}\left( \frac{1}{2}\left( \binom{%
n}{\left( n+1\right) /2}+\right. \right. \right.  \\
\left. \left. \left. \left. \binom{n}{\left( n-1\right) /2}2^{-\left(
n+1\right) /2}\left( n-1\right) /2-\binom{n}{\left( n+3\right) /2}2^{-\left(
n+3\right) /2}\left( n+3\right) /2\right) \right) \right) \right) .
\end{multline}%
Given that a function is equally likely to be picked from $M_{0,1}^{1}\left(
n,r,z,v\right) $ as from $M_{0,2}^{1}\left( n,r,z,v\right) ,$ the expected
average sensitivity is the average of equations (\ref{eqsf1}) and (\ref%
{eqsf2}).

\end{proof}

\begin{acknowledgement}
I am grateful to the Institute for Systems Biology (ISB) for their support.
\end{acknowledgement}

\bibliographystyle{plain}
\bibliography{refs}

\begin{thebibliography}{10}

\bibitem{aracena3}
Julio Aracena, Jacques Demongeot, and Eric Goles.
\newblock On limit cycles of monotone functions with symmetric connection
  graph.
\newblock {\em Theoretical Computer Science}, 322(2):237--244, 2004.

\bibitem{aracena1}
Julio Aracena, Maximilien Gadouleau, Adrien Richard, and Lilian Salinas.
\newblock Fixing monotone boolean networks asynchronously.
\newblock {\em Information and Computation}, page 104540, 2020.

\bibitem{aracena2}
Julio Aracena, Adrien Richard, and Lilian Salinas.
\newblock Number of fixed points and disjoint cycles in monotone boolean
  networks.
\newblock {\em SIAM Journal on Discrete Mathematics}, 31(3):1702--1725, 2017.

\bibitem{bernasconi2}
Anna Bernasconi.
\newblock Sensitivity vs. block sensitivity (an average-case study).
\newblock {\em Information processing letters}, 59(3):151--157, 1996.

\bibitem{bernasconi}
Anna Bernasconi, Carsten Damm, and Igor Shparlinski.
\newblock The average sensitivity of square-freeness.
\newblock {\em computational complexity}, 9(1):39--51, 2000.

\bibitem{boppana}
Ravi~B Boppana.
\newblock The average sensitivity of bounded-depth circuits.
\newblock {\em Information processing letters}, 63(5):257--261, 1997.

\bibitem{bshouty}
Nader~H Bshouty and Christino Tamon.
\newblock On the fourier spectrum of monotone functions.
\newblock {\em Journal of the ACM (JACM)}, 43(4):747--770, 1996.

\bibitem{derridapomeau}
Bernard Derrida and Yves Pomeau.
\newblock Random networks of automata: a simple annealed approximation.
\newblock {\em EPL (Europhysics Letters)}, 1(2):45, 1986.

\bibitem{friedgut}
Ehud Friedgut.
\newblock Boolean functions with low average sensitivity depend on few
  coordinates.
\newblock {\em Combinatorica}, 18(1):27--35, 1998.

\bibitem{friedgutkalai}
Ehud Friedgut and Gil Kalai.
\newblock Every monotone graph property has a sharp threshold.
\newblock {\em Proceedings of the American mathematical Society},
  124(10):2993--3002, 1996.

\bibitem{harris}
Stephen~E Harris, Bruce~K Sawhill, Andrew Wuensche, and Stuart Kauffman.
\newblock A model of transcriptional regulatory networks based on biases in the
  observed regulation rules.
\newblock {\em Complexity}, 7(4):23--40, 2002.

\bibitem{kahn}
J.~Kahn, G.~Kalai, and N.~Linial.
\newblock The influence of variables on boolean functions.
\newblock In {\em Proceedings of the 29th Annual Symposium on Foundations of
  Computer Science}, SFCS ’88, page 68–80, USA, 1988. IEEE Computer
  Society.

\bibitem{kauffman1969a}
Stuart~A Kauffman.
\newblock Metabolic stability and epigenesis in randomly constructed genetic
  nets.
\newblock {\em Journal of theoretical biology}, 22(3):437--467, 1969.

\bibitem{kauffmanbook}
Stuart~A Kauffman.
\newblock {\em The origins of order: Self-organization and selection in
  evolution}.
\newblock Oxford University Press, USA, 1993.

\bibitem{korshunov}
A.D. Korshunov.
\newblock On the number of monotone boolean functions.
\newblock {\em Problemy Kibernetiki}, 38:5--108, 1981.

\bibitem{korshunov2003monotone}
Aleksei~Dmitrievich Korshunov.
\newblock Monotone boolean functions.
\newblock {\em Russian Mathematical Surveys}, 58(5):929--1001, 2003.

\bibitem{ilya-special}
Aleksei~Dmitrievich Korshunov and Ilya Shmulevich.
\newblock The number of special monotone boolean functions and statistical
  properties of stack filters.
\newblock {\em Diskretnyi Analiz i Issledovanie Operatsii}, 7(3):17--44, 2000.

\bibitem{ilya-normal}
Aleksey~D Korshunov and Ilya Shmulevich.
\newblock On the distribution of the number of monotone boolean functions
  relative to the number of lower units.
\newblock {\em Discrete mathematics}, 257(2-3):463--479, 2002.

\bibitem{li}
Yuan Li, John~O Adeyeye, David Murrugarra, Boris Aguilar, and Reinhard
  Laubenbacher.
\newblock Boolean nested canalizing functions: A comprehensive analysis.
\newblock {\em Theoretical Computer Science}, 481:24--36, 2013.

\bibitem{luque}
Bartolo Luque and Ricard~V Sol{\'e}.
\newblock Lyapunov exponents in random boolean networks.
\newblock {\em Physica A: Statistical Mechanics and its Applications},
  284(1-4):33--45, 2000.

\bibitem{melliti}
Tarek Melliti, Damien Regnault, Adrien Richard, and Sylvain Sen{\'e}.
\newblock Asynchronous simulation of boolean networks by monotone boolean
  networks.
\newblock In {\em International Conference on Cellular Automata}, pages
  182--191. Springer, 2016.

\bibitem{mesot-teuscher}
Bertrand Mesot and Christof Teuscher.
\newblock Critical values in asynchronous random boolean networks.
\newblock In {\em European Conference on Artificial Life}, pages 367--376.
  Springer, 2003.

\bibitem{mossel}
Elchanan Mossel and Ryan O'Donnell.
\newblock On the noise sensitivity of monotone functions.
\newblock {\em Random Structures \& Algorithms}, 23(3):333--350, 2003.

\bibitem{munoz}
Miguel~A Munoz.
\newblock Colloquium: Criticality and dynamical scaling in living systems.
\newblock {\em Reviews of Modern Physics}, 90(3):031001, 2018.

\bibitem{serra}
Roberto Serra, Marco Villani, and Alessandro Semeria.
\newblock Genetic network models and statistical properties of gene expression
  data in knock-out experiments.
\newblock {\em Journal of theoretical biology}, 227(1):149--157, 2004.

\bibitem{shi}
Yaoyun Shi.
\newblock Lower bounds of quantum black-box complexity and degree of
  approximating polynomials by influence of boolean variables.
\newblock {\em Information Processing Letters}, 75(1-2):79--83, 2000.

\bibitem{ilyansip95}
I~Shmulevich, TM~Sellke, M~Gabbouj, and EJ~Coyle.
\newblock Stack filters and free distributive lattices.
\newblock In {\em Proceedings of the 1995 IEEE workshop on nonlinear signal and
  image processing}, pages 829--930, 1995.

\bibitem{ilya-arxiv-MBF}
Ilya Shmulevich.
\newblock Average sensitivity of typical monotone boolean functions.
\newblock {\em arXiv preprint math/0507030}, 2005.

\bibitem{pbn1}
Ilya Shmulevich, Edward~R Dougherty, Seungchan Kim, and Wei Zhang.
\newblock Probabilistic boolean networks: a rule-based uncertainty model for
  gene regulatory networks.
\newblock {\em Bioinformatics}, 18(2):261--274, 2002.

\bibitem{ilya-kauffman}
Ilya Shmulevich and Stuart~A Kauffman.
\newblock Activities and sensitivities in boolean network models.
\newblock {\em Physical review letters}, 93(4):048701, 2004.

\bibitem{post}
Ilya Shmulevich, Harri L{\"a}hdesm{\"a}ki, Edward~R Dougherty, Jaakko Astola,
  and Wei Zhang.
\newblock The role of certain post classes in boolean network models of genetic
  networks.
\newblock {\em Proceedings of the National Academy of Sciences},
  100(19):10734--10739, 2003.

\bibitem{ilya-almostall}
Ilya Shmulevich, Olli Yli-Harja, Jaakko Astola, and Aleksey Korshunov.
\newblock On the robustness of the class of stack filters.
\newblock {\em IEEE Transactions on Signal Processing}, 50(7):1640--1649, 2002.

\bibitem{zhang}
Shengyu Zhang.
\newblock Note on the average sensitivity of monotone boolean functions.
\newblock {\em Preprint}, page~4, 2011.

\end{thebibliography}

\end{document}